\documentclass[11pt,titlepage]{article}
\usepackage{amssymb}
\usepackage{latexsym}
\usepackage[all]{xy}
\usepackage{array}
\newtheorem{defi}{Definition}[section]
\newtheorem{lm}{Lemma}[section]
\newtheorem{thm}{Theorem}[section]
\newtheorem{prop}{Proposition}[section]
\newcommand{\bew}[1]{\noindent{\bf Proof:} #1}
\newcommand{\opm}[1]{\noindent{\bf Remark:} #1}
\newcommand{\vb}[1]{\noindent{\bf Examples:}\begin{enumerate} #1 \end{enumerate}}
\newcommand{\subtilde}[1]{\oalign{$#1$\crcr\hidewidth\vbox to .2ex{\hbox{\char'176}\vss}\hidewidth}}
\newcommand{\colimit}{\oalign{$\lim$\crcr\hidewidth\vbox to .2ex{\hbox{$\rightarrow$}\vss}\hidewidth}}

\makeatletter
\renewcommand\section{\@startsection {section}{1}{\z@}%
                                   {-3ex \@plus -0ex \@minus
                                   -1ex}
                                   {2.3ex \@plus.2ex}%
                                   {\normalfont\LARGE\bfseries 
Lecture\ }}
\makeatother

\title{Introduction to the Language of Stacks and Gerbes}
\author{I. Moerdijk, University of Utrecht}
\date{}

\begin{document}

\maketitle

{\normalfont\LARGE\bfseries Introduction}

\vspace{0.4cm}
\noindent This is the text of a series of four one hour lectures given as part of
the ``Third Lisbon Summer Lectures in Geometry'', which took place in the
summer of 2002 at the Instituto Tecnico Superior in Lisbon. The
lectures were aimed at an audience consisting of students with some
background in topology (but not necessarily in algebraic geometry). The
purpose of the lectures was to give a quick introduction to gerbes,
with an emphasis on non-abelian \v{C}ech cohomology.

The lectures start with the concept of a sheaf on a space, and the
well-known description of principal bundles by non-abelian cocycles in
degree 1. The second lecture introduces the three notions of fibered
category, prestack and stack, with some emphasis on the analogy with
the three concepts of presheaf, separated presheaf and sheaf. In the
third lecture, gerbes and bands (``liens'') are introduced, and the
description of equivalence classes of gerbes with a given band by
non-abelian cocycles of degree 2 is described in detail. Along the way,
it is shown that every equivalence class of such gerbes can be
represented by a sheaf of groupoids. In the fourth and last lecture,
some relations with differential geometry are discussed. In particular,
we describe the notion of bundle gerbe and its relation to extensions
of Lie groupoids, and we relate the degree 2 cohomology class of the
previous lecture to De Rham cohomology in the abelian case.

Given the fact that the material had to fit into four hours of
lectures, the presentation is necessarily rather concise. Nonetheless,
I hope that it represents a rather efficient introduction to some of
the main concepts, which could be of use to some students.

Apart from perhaps an occasional remark, there is nothing original in
these notes, and much of it can already be found in the classical
reference Giraud {\bf [Gi]} for non-abelian cohomology. This applies in
particular to the material of first two lectures. The presentation in
the third lecture is a variation on the description in Breen {\bf [Br]}. The De
Rham theory of bundle gerbes is based on Murray {\bf [Mu]}, while the definition
of general bundle gerbes in terms of extensions of Lie groupoids is
based on Moerdijk {\bf [Mo]}.

\newpage

\section{Sheaves and Torsors}

\subsection*{Sheaves}
We assume familiarity with the basic categorical terminology of categories, functors and natural transformations.
Let $X$ be a topological space. The poset $\mathcal{O}(X)$ of open subsets of $X$ can be viewed as a category, with inclusions $i:U\hookrightarrow V$ as arrows (for $U\subseteq V$). We shall work with this category, although everything applies in fact to an arbitrary site.

A \emph{presheaf} (of sets) is a contravariant functor $P$ from $\mathcal{O}(X)$ into Sets,
\begin{displaymath}
P:\quad \mathcal{O}(X)^{op}\rightarrow \mathrm{Sets}.
\end{displaymath}
Thus, we have a set $P(U)$ for each open $U\subseteq X$, and a ``restriction operation'' $P(i):P(V)\rightarrow P(U)$ whenever $i:U\hookrightarrow V$ is an inclusion. For $a\in P(V)$ we usually write $i^*(a)$ or $a|U$ for $P(i)(a)$. A morphism between presheaves is a natural transformation. This defines the category of presheaves on $X$.

Let $U=\bigcup U_{\alpha}$ be an open cover. A \emph{compatible family} in the presheaf $P$ consists of elements $a_{\alpha}\in P(U_{\alpha})$ which agree on overlaps, $a_{\alpha}|U_{\alpha\beta}=a_{\beta}|U_{\alpha\beta}$ for any two indices $\alpha$ and $\beta$ (here $U_{\alpha\beta}=U_{\alpha}\cap U_{\beta}$). A presheaf $P$ is called a \emph{sheaf} if for any such cover $U=\bigcup U_{\alpha}$, and for any compatible family $\{a_{\alpha}\}$, there is a \emph{unique} ``amalgamate'' $a\in P(U)$ for which $a|U_{\alpha}=a_{\alpha}$ (for each $\alpha$). If there is \emph{at most} one such $a$, the presheaf is called \emph{separated}. With the same natural transformations as morphisms, this defines two smaller categories, that of sheaves on $X$ and of separated presheaves on $X$.
In a similar way, one defines (pre-)sheaves of groups, rings, etc.

\vspace{0.4cm}
\vb{
\item Let $p:Y\rightarrow X$ be a map. Define $\mathcal{S}(p)(U)=\{f:U\rightarrow Y|\ p\circ f=\mathrm{identity,\ }f\mathrm{\ continuous}\}$. Then $\mathcal{S}(p)$ is a sheaf, the \emph{sheaf of sections} of $p$.\label{shsections}
\item If $G$ is a Lie group and $X$ is a manifold, the sheaf of smooth sections of the projection $G\times X\rightarrow X$ is a sheaf of groups, denoted $\subtilde{G}_X$, or just $\subtilde{G}$.\label{lie}
}

\begin{thm}
For a given space $X$, the inclusion functors
\[
\mathrm{Sheaves}\hookrightarrow\mathrm{Separated\ presheaves}\hookrightarrow \mathrm{Presheaves\ on\ }X
\]
have left adjoints.
\end{thm}

For a presheaf $P$ on $X$ we will write $\bar P$ for the associated separated presheaf and $\mathbf{a}(P)$ for the \emph{associated sheaf}, the composition of the two adjoints.
\\

\bew{(sketch)
\begin{enumerate}
\item If $P$ is a presheaf, define $\bar P(U)=P(U)/\sim$ where $a\sim b$ iff there is a cover $U=\bigcup U_{\alpha}$ with $a|U_{\alpha}=b|U_{\alpha}$.
\item If $P$ is a \emph{separated presheaf}, define $\mathbf{a}(P)(U)$ to be the set of equivalence classes of compatible families for covers $U=\bigcup U_{\alpha}$, the equivalence relation being given by refinement of covers.
\end{enumerate}
}

\vspace{0.4cm}
\opm{For a separated presheaf $P$ its associated sheaf is characterized uniquely up to isomorphism as a sheaf $Q$ equipped with an embedding $\eta:P\hookrightarrow Q$ which is \emph{locally surjective} (i.e. for any $b\in Q(U)$ there is a cover $U=\bigcup U_{\alpha}$ such that each $b|U_{\alpha}$ is in the image of $\eta:P(U_{\alpha})\hookrightarrow Q(U_{\alpha})$).
}

\subsection*{Etale spaces}
An \'etale space over $X$ is a space $E$ equipped with a continuous map $\phi:E\rightarrow X$ which is a \emph{local homeomorphism} (i.e. for any $y\in E$ there are open neighbourhoods $V_y\subseteq E$ of $y$ and $U_x$ of $x=\phi(y)$ in $X$ for which $\phi$ restricts to a homeomorphism $V_y\stackrel{\sim}{\rightarrow}U_x$). For two such \'etale spaces $\phi:E\rightarrow X$ and $\psi:F\rightarrow X$ we consider continuous maps $f:E\rightarrow F$ for which $\psi f=\phi$. This defines a category of \'etale spaces over $X$. There is an obvious functor 
\begin{equation}
\mathcal{S}:\mathrm{Etale\ spaces}\rightarrow\mathrm{Sheaves}\label{siers}
\end{equation}
sending $\phi:E\rightarrow X$ to the sheaf of sections $\mathcal{S}(\phi)$.

\begin{thm}
This functor is an equivalence of categories.
\end{thm}

\bew{We sketch the construction of its inverse. Given a sheaf $P$ on $X$ and a point $x\in X$, write 
\[
P_x=\colimit_{x\in U} P(U)
\]
for the stalk of $P$ at $x$. For any open $U\ni x$, we write germ$_x: P(U)\rightarrow P_x$ for the canonical map. (In Example \ref{shsections} above $\mathcal{S}(p)_x$ is the set of germs at $x$ of sections of $p$.) Let $E(P)=\coprod_{x\in X}P_x$ be the set of all germs, with evident projection
\[
\pi:E(P)\rightarrow X.
\]
The topology on $E(P)$ is defined by basic opens
\[
B(a)=\{\mathrm{germ}_x(a)|x\in U\}
\]
where $U\subseteq X$ is open and $a\in P(U)$. This topology makes $E(P)$ into an \'etale space over $X$, and provides a functor which is inverse - up to natural isomorphisms - to (\ref{siers}).
}
\\

On the basis of this theorem we can identify sheaves and \'etale spaces over $X$, and we will often (tacitly) move from one context to another.

\vspace{0.4cm}
\noindent{\bf Remarks:}
\begin{enumerate}
\item For a sheaf $P$, elements of $P(U)$ correspond to sections of $E(P)\rightarrow X$ over $U$. For this reason one refers to elements of $P(X)$ as \emph{global sections} of $P$, and also write $\Gamma(P)$ or $\Gamma(X,P)$ for $P(X)$.
\item The product of sheaves is constructed as $(P\times Q)(U)=P(U)\times Q(U)$. The corresponding operation on \'etale spaces $E\rightarrow X$ and $F\rightarrow X$ is the fibered product $E\times_X F$ over $X$. In particular, we have $E(P)\times_X E(Q)\cong E(P\times Q)$.
\item A sheaf of groups on $X$ is ``the same'' as an \'etale space $G\rightarrow X$, equipped with a continuously varying group structure on each of its fibers $G_x$ given by unit and multiplication maps $u:X\rightarrow G$ and $m:G\times G\rightarrow G$ (maps \emph{over} $X$).\label{fibergroup}
\end{enumerate}

\subsection*{Torsors}
Let $\mathcal{G}$ be a sheaf of groups on $X$, and let $\mathcal{S}$ be a sheaf on $X$. An \emph{action} of $\mathcal{G}$ on $\mathcal{S}$ is a map of sheaves
\[
\mu:\mathcal{G}\times\mathcal{S}\rightarrow\mathcal{S}
\]
(with components $\mu_U:\mathcal{G}(U)\times\mathcal{S}(U)\rightarrow\mathcal{S}(U)$ denoted $\mu_U(g,a)=g\cdot a$) which satisfies the usual conditions for a (left) action. $\mathcal{S}$ is called a \emph{$\mathcal{G}$-torsor} if
\begin{enumerate}
\item $X=\bigcup\{U|\mathcal{S}(U)\neq\emptyset\}$\label{torconX}
\item For each open $U\subseteq X$ the action of $\mathcal{G}(U)$ on $\mathcal{S}(U)$ is free and transitive.\label{torconft}
\end{enumerate}

A \emph{morphism} $\mathcal{S}\rightarrow\mathcal{S}'$ between $\mathcal{G}$-torsors is a morphism of sheaves which commutes with the action. Any such morphism is automatically an \emph{isomorphism}.

This definition translates into \'etale spaces as follows: for an \'etale space $G\rightarrow X$ with fiberwise group structure as in Example \ref{fibergroup} above, a $G$-torsor is an \'etale space $E\rightarrow X$ equipped with an action $\mu:G\times_X E\rightarrow E$ for which
\begin{enumerate}
\item $E\rightarrow X$ is surjective
\item the map $(\mu,\pi_2):G\times_X E\rightarrow E\times_X E$ is a homeomorphism.
\end{enumerate}

\vb{
\item Let $L$ be a Lie group and let $\pi:P\rightarrow X$ be a principal bundle over a manifold $X$. Then the sheaf $\mathcal{S}(\pi)$ of smooth sections is an $\subtilde{L}$-torsor (cf. Example \ref{lie} on page \pageref{lie}).
\item For a sheaf $\mathcal{G}$ of groups, the group multiplication $\mathcal{G}\times\mathcal{G}\rightarrow\mathcal{G}$ makes $\mathcal{G}$ itself into a $\mathcal{G}$-torsor. A $\mathcal{G}$-torsor $\mathcal{S}$ is isomorphic to this ``trivial'' $\mathcal{G}$-torsor iff $\mathcal{S}(X)\neq\emptyset$, i.e. iff $\mathcal{S}$ has a global section.
}

\subsection*{Non-abelian cohomology in degree 1}
Let $\mathcal{G}$ be a sheaf of groups on $X$, and let $\mathcal{U}=\{U_{\alpha}\}$ be an open cover of $X$. A \emph{1-cocycle} with values in $\mathcal{G}$ is a family $\underline{g}=\{g_{\alpha\beta}\},\ g_{\alpha\beta}\in\mathcal{G}(U_{\alpha\beta})$,  satisfying the cocycle indentity
\begin{equation}
g_{\alpha\beta}g_{\beta\gamma}=g_{\alpha\gamma}\quad\mathrm{on}\quad U_{\alpha\beta\gamma}\label{coid}
\end{equation}
for any three indices $\alpha,\beta,\gamma$. (More precisely we should have written \linebreak $(g_{\alpha\beta}|U_{\alpha\beta\gamma})\cdot(g_{\beta\gamma}|U_{\alpha\beta\gamma})=(g_{\alpha\gamma}|U_{\alpha\beta\gamma})$, or even $\mu_{U_{\alpha\beta\gamma}}(g_{\alpha\beta}|U_{\alpha\beta\gamma},g_{\beta\gamma}|U_{\alpha\beta\gamma})=g_{\alpha\gamma}|U_{\alpha\beta\gamma}$, but we usually just simply abuse notation as in (\ref{coid}).) Two such cocycles $\underline{g}$ and $\underline{h}$ are \emph{equivalent} if there are $f_{\alpha}\in\mathcal{G}(U_{\alpha})$ for which 
\[
f_{\alpha}g_{\alpha\beta}=h_{\alpha\beta}f_{\beta}.
\]
The set of equivalence classes is denoted $H^1(\mathcal{U},\mathcal{G})$, and the equivalence class of $\underline{g}$ is denoted $[\underline{g}]$. The colimit (over refinement of covers) is denoted
\[
\check{H}^1(X,\mathcal{G}):=\colimit_{\mathcal{U}}H^1(\mathcal{U},\mathcal{G}).
\]

\subsection*{Cocycle description of torsors}
Let $\mathcal{S}$ be a $\mathcal{G}$-torsor over $X$. By condition \ref{torconX} we can choose a cover $X=\bigcup U_{\alpha}$ and sections $s_{\alpha}\in\mathcal{S}(U_{\alpha})$. By condition \ref{torconft} there are unique $g_{\alpha\beta}\in\mathcal{G}(U_{\alpha\beta})$ for which
\[
g_{\alpha\beta}\cdot s_{\beta}=s_{\alpha}\quad\mathrm{on}\quad U_{\alpha\beta}
\]
(more precisely $g_{\alpha\beta}\cdot (s_{\beta}|U_{\alpha\beta})=(s_{\alpha}|U_{\alpha\beta})\in\mathcal{S}(U_{\alpha\beta})$). Then $\underline{g}=\{g_{\alpha\beta}\}$ is a cocycle. Suppose we had chosen different sections $t_{\alpha}\in\mathcal{S}(U_{\alpha})$, giving a cocycle $\underline{h}=\{h_{\alpha\beta}\},\ h_{\alpha\beta}\cdot t_{\beta}=t_{\alpha}\ \mathrm{on}\ U_{\alpha\beta}$. Let $f_{\alpha}\in\mathcal{G}(U_{\alpha})$ be the unique element with $f_{\alpha}\cdot s_{\alpha}=t_{\alpha}$. Then $f_{\alpha}\cdot g_{\alpha\beta}=h_{\alpha\beta}\cdot f_{\beta}$, showing that $[\underline{g}]=[\underline{h}]$. If we also allow the cover $\{U_{\alpha}\}$ to vary, we obtain a welldefined class $[\mathcal{S}]\in\check{H}^1(X,\mathcal{G})$.

\subsection*{Torsors from cocycles}
Given a cocycle $\{g_{\alpha\beta}\}$ as before, we can construct a torsor. For variation we give a construction in terms of \'etale spaces. Let $\pi:G\rightarrow X$ be the \'etale space corresponding to the sheaf $\mathcal{G}$, so that $g_{\alpha\beta}$ is a section of $\pi$ over $U_{\alpha\beta}$. Write $G_{U_{\alpha}}=\pi^{-1}(U_{\alpha})$, and let
\[
E=\left(\coprod_{\alpha}G_{U_{\alpha}}\right)/\sim.
\]
Writing points of $\coprod_{\alpha}G_{U_{\alpha}}$ (somewhat redundantly) as triples $(x,\alpha,g)$ where $x\in U_{\alpha}$ and $g\in G_x=\pi^{-1}(x)$, the equivalence relation $\sim$ here is defined by
\[
(x,\alpha,g)\sim(x,\beta,gg_{\alpha\beta}(x)).
\]
Then $E$ with its obvious left $G$-action is a $G$-torsor.
The previous two constructions are mutually inverse and show:
\begin{thm}
There is a bijective correspondence between isomorphism \linebreak classes of $\mathcal{G}$-torsors and cohomology classes in $\check{H}^1(X,\mathcal{G})$.
\end{thm}

\newpage

\section{Stacks}
Again we work over a topological space $X$ and its category of open sets and inclusions $\mathcal{O}(X)$. Also, to understand the notation, recall that if
\[
\xymatrix{\mathbb{A}\ar[r]^{F}&\mathbb{B} \ar@<0.5ex>[r]^{G} \ar@<-0.5ex>[r]_{H}&\mathbb{C} \ar[r]^{K}&\mathbb{D}}
\]
are functors, then any natural transformation $\tau:G\rightarrow H$ induces natural transformations $K\tau:KG\rightarrow KH$ and $\tau_F$ or $\tau F:GF\rightarrow HF$.
\begin{defi}
A \emph{fibered category} $\mathbb{F}$ over $X$ consists of
\begin{itemize}
\item a category $\mathbb{F}(U)$ for each open $U\subseteq X$,
\item a functor $i^*:\mathbb{F}(U)\rightarrow\mathbb{F}(V)$ for each inclusion $i:V\hookrightarrow U$ in $\mathcal{O}(X)$,
\item a natural isomorphism
\[
\tau=\tau_{i,j}:(ij)^*\rightarrow j^*i^*
\]
for each pair of inclusions $W\stackrel{j}{\hookrightarrow}V\stackrel{i}{\hookrightarrow}U$.
\end{itemize}
Moreover, for any three composable inclusions $N\stackrel{k}{\hookrightarrow}W\stackrel{j}{\hookrightarrow}V\stackrel{i}{\hookrightarrow}U$, the diagram 
\[
\xymatrix{(ijk)^* \ar[r]^{\tau_{ij,k}} \ar[d]_{\tau_{i,jk}} & k^*(ij)^* \ar[d]^{k^*\tau_{i,j}} \\ (jk)^*i^* \ar[r]^{\tau_{j,k}i^*} &k^*j^*i^*} 
\]
should commute.
\end{defi}

Sometimes we write more explicitly $(\mathbb{F},\tau)$ for the fibered category $\mathbb{F}$. Sometimes we write $a|V$ for $i^*(a)$, for $i:V\hookrightarrow U$ and $a\in\mathbb{F}(U)$.

\vspace{0.4cm}
\noindent{\bf Example:}
A presheaf of categories is the same thing as a fibered category $(\mathbb{F},\tau)$ in which all $\tau_{i,j}$ are identity transformations. In a way, a fibered category is a ``presheaf up to isomorphisms''.

\begin{defi}
Let $\mathbb{F}$ and $\mathbb{G}$ be fibered categories over $X$. A morphism $\phi:\mathbb{F}\rightarrow\mathbb{G}$ of fibered categories consists of
\begin{itemize}
\item a functor $\phi(U):\mathbb{F}(U)\rightarrow\mathbb{G}(U)$ for each open $U\subseteq X$,
\item a natural isomorphism
\[
\alpha_i:\phi_V i^*\stackrel{\sim}{\rightarrow}i^*\phi_U
\]
for each inclusion $i:V\hookrightarrow U$,
\end{itemize}
and these should satisfy a compatibility condition with respect to the $\tau$'s: for inclusions $W\stackrel{j}{\hookrightarrow}V\stackrel{i}{\hookrightarrow}U$, the diagram
\[
\xymatrix{& \phi_W(ij)^* \ar[dl]_{\phi_W\tau} \ar[rr]^{\alpha_{ij}} & & (ij)^*\phi_U \ar[dr]^{\tau\phi_U} \\\phi_Wj^*i^* \ar[rr]^{(\alpha_j)i^*} & &j^*\phi_Vi^* \ar[rr]^{j^*\alpha_i} & & j^*i^*\phi_U
}
\]
should commute.
\end{defi}

\begin{defi}
Such a morphism $\phi$ (or $(\phi,\alpha)$, more precisely) is called a \emph{strong equivalence} if every $\phi_U$ is an equivalence of categories. It is called a \emph{weak equivalence} if every $\phi_U$ is fully faithful, and ``locally surjective'' on objects, in the sense that for every object $a\in\mathbb{G}(U)$ and every $x\in U$, there exists a neighbourhood $V$, $x\in V\stackrel{i}{\hookrightarrow}U$, and an object $b\in\mathbb{F}(V)$ such that $\phi_V(b)$ is isomorphic in $\mathbb{G}(V)$ to $i^*(a)$.
\end{defi}

\opm{
(To be skipped on first reading) Suppose $\phi=(\phi,\alpha)$ and $\psi=(\psi,\beta):\mathbb{F}\rightarrow\mathbb{G}$ are two morphisms between fibered categories. A \emph{fibered transformation} $\mu:\phi\rightarrow\psi$ consists of (ordinary) natural transformations
\[
\mu_U:\phi_U\rightarrow\psi_U,
\]
one for each open $U\subseteq X$. These $\mu$'s are required to be compatible with the $\alpha$'s and the $\beta$'s in the following sense: for any inclusion $i:V\hookrightarrow U$, the diagram of natural transformations
\[
\xymatrix{\phi_Vi^* \ar[r]^{\alpha_i} \ar[d]_{\mu_Vi^*} & i^*\phi_U \ar[d]^{i^*\mu_U} \\ \psi_Vi^* \ar[r]^{\beta_i} &i^*\psi_U} 
\]
should commute. We say $\mu$ is a \emph{fibered isomorphism} if each $\mu_U$ is a natural isomorphism.
\\
Fibered categories, with their morphisms and fibered transformations, form a 2-category. (In the sequel we will suppress 2-categorical details.)
}

\begin{lm}\label{hom}
\begin{enumerate}
\item Let $\mathbb{F}$ be a fibered category over $X$ and let $a,b$ be objects in $\mathbb{F}(U)$. Then the assignment
\[
V\mapsto\mathrm{Hom}_{\mathbb{F}(V)}(i^*a,i^*b)
\]
for $i:V\hookrightarrow U$, defines a presheaf on $U$. This presheaf is denoted $\mathrm{\underline{Hom}}_{\mathbb{F}}(a,b)$.
\item Any morphism $\phi:\mathbb{F}\rightarrow\mathbb{G}$ of fibered categories induces a morphism of presheaves on $U$:
\[
\phi_{a,b}:\mathrm{\underline{Hom}}_{\mathbb{F}}(a,b)\rightarrow\mathrm{\underline{Hom}}_{\mathbb{G}}(\phi_U(a),\phi_U(b)).
\]
\end{enumerate}
\end{lm}

\bew{
With the following explicit definitions, this is a matter of straightforward checking:\\
(i) The restriction maps of the presheaf are given, for $W\stackrel{j}{\hookrightarrow}V\stackrel{i}{\hookrightarrow}U$, by the composition
\[
\mathrm{Hom}_{\mathbb{F}(V)}(i^*a,i^*b)\stackrel{j^*}{\rightarrow}\mathrm{Hom}_{\mathbb{F}(W)}(j^*i^*a,j^*i^*b)\stackrel{\tau^*}{\rightarrow}\mathrm{Hom}_{\mathbb{F}(W)}((ij)^*a,(ij)^*b)
\]
where $\tau^*=\tau_{i,j}^*$ is ``conjugation by $\tau$'', mapping an arrow $f$ to $(\tau_{i,j})_b^{-1}f(\tau_{i,j})_a$.\\
(ii) The component of $\phi_{a,b}$ at $i:V\hookrightarrow U$ is the composition
\[
\mathrm{Hom}_{\mathbb{F}(V)}(i^*a,i^*b)\stackrel{\phi_V}{\rightarrow}\mathrm{Hom}_{\mathbb{G}(V)}(\phi_Vi^*a,\phi_Vi^*b)\stackrel{(\alpha_i)_*}{\rightarrow}\mathrm{Hom}_{\mathbb{G}(V)}(i^*\phi_Ua,i^*\phi_Ub)
\]
where $(\alpha_i)_*$ is conjugation by $\alpha_i,\ g\mapsto(\alpha_i)_bg(\alpha_i)_a^{-1}$.
}

\begin{defi}
A fibered category $\mathbb{F}$ over $X$ is called a \emph{prestack} if, for any objects $a,b\in\mathbb{F}(U)$, the presheaf $\mathrm{\underline{Hom}}_{\mathbb{F}}(a,b)$ is a sheaf.
\end{defi}

\noindent{\bf Example:}
Let $P$ be a presheaf of sets on $U$, viewed as a presheaf of categories with identity arrows only, hence viewed as a special kind of fibered category over $U$. Then $P$ is a prestack iff it is a separated presheaf.

\begin{prop}
For any space $X$ the inclusion
\[
\mathrm{Prestacks}\hookrightarrow\mathrm{Fibered\ categories\ over\ }X
\]
has a left adjoint, denoted: $\mathbb{F}\mapsto\bar{\mathbb{F}}$.
\end{prop}

\bew{
This is easy, using the lemma: given a fibered category $\mathbb{F}$ over $X$, one defines a new fibered category $\bar{\mathbb{F}}$ by letting $\bar{\mathbb{F}}(U)$ have the same objects as $\mathbb{F}(U)$ and arrows defined by
\[
\mathrm{Hom}_{\bar{\mathbb{F}}}(a,b)=\Gamma{\bf a}(\mathrm{\underline{Hom}}_{\mathbb{F}}(a,b)),
\]
the set of global sections of the sheaf associated to the presheaf $\mathrm{\underline{Hom}}_{\mathbb{F}}(a,b)$ on $U$. Further details are straightfoward.
}

\begin{defi}
Let $\mathbb{F}$ be a fibered category over $X$ and let $\mathcal{U}=\{U_{\alpha}\}_{\alpha\in A}$ be an open cover of an open set $U\subseteq X$. The category $\mathrm{Des}(\mathcal{U},\mathbb{F})$ of ``descent data'' has
\begin{itemize}
\item as \emph{objects}:\quad systems $(a,\theta)=(\{a_{\alpha}\},\{\theta_{\alpha\beta}\})$, where each $a_{\alpha}$ is an object of $\mathbb{F}(U_{\alpha})$ and
\[
\theta_{\alpha\beta}:a_{\beta}|U_{\alpha\beta}\stackrel{\sim}{\rightarrow}a_{\alpha}|U_{\alpha\beta}
\]
is an isomorphism (here $a_{\beta}|U_{\alpha\beta}=i^*(a_{\beta})$ for the inclusion $i:U_{\alpha\beta}\hookrightarrow U_{\beta}$; etc). These isomorphisms are required to satisfy the \emph{cocycle condition}
\[
\theta_{\alpha\alpha}=1,\quad\theta_{\alpha\beta}\circ\theta_{\beta\gamma}=\theta_{\alpha\gamma}\ (\mathrm{in\ }\mathbb{F}(U_{\alpha\beta\gamma})).
\]
\item as \emph{arrows}:\quad $(a,\theta)\stackrel{f}{\rightarrow}(b,\rho)$, families of arrows $f_{\alpha}:a_{\alpha}\rightarrow b_{\alpha}$ in $\mathbb{F}(U_{\alpha})$ for which $\rho_{\alpha\beta}f_{\beta}=f_{\alpha}\theta_{\alpha\beta}$, i.e. each diagram
\[
\xymatrix {a_{\beta}|U_{\alpha\beta}\ar[d]_{\theta_{\alpha\beta}}\ar[r]^{f_{\beta}} & b_{\beta}|U_{\alpha\beta}\ar[d]^{\rho_{\alpha\beta}} \\ a_{\alpha}|U_{\alpha\beta}\ar[r]^{f_{\alpha}} & b_{\beta}|U_{\alpha\beta}
}
\]
commutes.
\end{itemize}
\end{defi}

\opm{
There is an obvious functor
\[
D=D(\mathcal{U},\mathbb{F}):\mathbb{F}(U)\rightarrow\mathrm{Des}(\mathcal{U},\mathbb{F});
\]
This functor is fully faithful (for each open $U$ and each cover $U=\bigcup U_{\alpha}$) iff $\mathbb{F}$ is a prestack.
}

\begin{defi}
The fibered category $\mathbb{F}$ is called a \emph{stack} if each such functor $D:\mathbb{F}(U)\rightarrow\mathrm{Des}(\mathcal{U},\mathbb{F})$ is an equivalence of categories.
\end{defi}

\vspace{0.4 cm}
\noindent{\bf Remarks:}
\begin{enumerate}
\item A fibered category $\mathbb{F}$ is a stack iff it is a prestack and, for every cover $U=\bigcup U_{\alpha}$, any object $(a,\theta)$ in $\mathrm{Des}(\mathcal{U},\mathbb{F})$ is isomorphic to an object of the form $D(b)$ for some $b\in\mathbb{F}(U)$ (where $D:\mathbb{F}(U)\rightarrow\mathrm{Des}(\mathcal{U},\mathbb{F})$ and $\mathcal{U}=\{U_{\alpha}\}$).
\item The construction of the category $\mathrm{Des}(\mathcal{U},\mathbb{F})$ is functorial: a morphism $\phi:\mathbb{F}\rightarrow\mathbb{G}$ of fibered categories induces a functor $\mathrm{Des}(\mathcal{U},\mathbb{F})\rightarrow\mathrm{Des}(\mathcal{U},\mathbb{G})$.
\item For a cover $\mathcal{U}=\{U_{\alpha}\}$ of $U$ and an inclusion $i:V\hookrightarrow U$, we obtain a cover $\{V_{\alpha}\}$ of $V$ and a canonical functor
\[
i^*:\mathrm{Des}(\mathcal{U},\mathbb{F})\rightarrow\mathrm{Des}(\mathcal{V},\mathbb{F}),
\]
by setting $V_{\alpha}=U_{\alpha}\cap V$ with inclusion $i_{\alpha}:V_{\alpha}\hookrightarrow U_{\alpha}$; the functor $i^*$ sends an object $(a,\theta)=\{a_{\alpha},\theta_{\alpha\beta}\}$ to the family of objects $i^*_{\alpha}(a_{\alpha})\in\mathbb{F}(V_{\alpha})$, together with the isomorphisms
\[
i^*_{\beta}(a_{\beta})|V_{\alpha\beta}\rightarrow i^*_{\alpha}(a_{\alpha})|V_{\alpha\beta}
\]
obtained from $\theta_{\alpha\beta}$ and the structure maps $\tau$ of $\mathbb{F}$. Note that for this functor $i^*$, the diagram
\[
\xymatrix{\mathbb{F}(U) \ar[d]_{i^*} \ar[r]^{D} & \mathrm{Des}(\mathcal{U},\mathbb{F}) \ar[d]^{i^*} \\ \mathbb{F}(V) \ar[r]^{D} & \mathrm{Des}(\mathcal{V},\mathbb{F})
}
\]
commutes.\label{rem3}
\end{enumerate}

\vspace{0.4cm}
\vb{
\item The assignment, to each open set $U\subseteq X$, of the category Vect$(U)$ of vector bundles over $U$ is a stack, as is the assignment of the category Sheaves$(U)$ of sheaves on $U$. If $\mathcal{G}$ is a sheaf of groups on $U$, one obtains a stack Tor$(\mathcal{G})$ by assigning to each $U$ the category Tor$(\mathcal{G}|U)$ of $\mathcal{G}|U$-torsors over $U$.
\item Related to the last example, one can also define a fibered category $U\mapsto\mathcal{G}(U)$, where the group $\mathcal{G}$ is viewed as a category with only one object. Observe that $\mathcal{G}$ is a prestack but not a stack. Observe that there is a natural weak equivalence $\mathcal{G}\rightarrow$Tor$(\mathcal{G}):$ for an open $U\subseteq X$, the unique object of $\mathcal{G}(U)$ is sent to the trivial $\mathcal{G}|U$-torsor over $U$. This means that Tor$(\mathcal{G})$ is the \emph{associated stack} of $\mathcal{G}$; see below.\label{asstack}
}

\vspace{0.4 cm}
\noindent{\bf Remarks:}
\begin{enumerate}
\item Let $w:\mathbb{F}\rightarrow\mathbb{G}$ be a weak equivalence of prestacks over $X$. If $\mathbb{F}$ is a stack then $w$ is also a strong equivalence.
\item Let $w:\mathbb{F}\rightarrow\mathbb{G}$ be a weak equivalence of prestacks. Then for any map $\psi:\mathbb{F}\rightarrow\mathbb{S}$ into a stack, there is an extension $\tilde{\psi}:\mathbb{G}\rightarrow\mathbb{S}$ with $\tilde{\psi}w\simeq\psi$ (natural isomorphism). Moreover, this extension $\tilde{\psi}$ is itself unique up to fibered isomorphism.
\item Let $\phi:\mathbb{F}\rightarrow\mathbb{G}$ be a strong equivalence. Then $\mathbb{F}$ is a (pre-) stack iff $\mathbb{G}$ is.
\end{enumerate}
The proofs of these assertions can safely be left as exercises.

\begin{defi}
Let $\mathbb{F}$ be a prestack. An associated stack for $\mathbb{F}$ consists of a stack $\hat{\mathbb{F}}$ together with a weak equivalence $\mathbb{F}\stackrel{w}{\rightarrow}\hat{\mathbb{F}}$.
\end{defi}

By the previous remarks, this associated stack is unique up to strong equivalence, and has a universal property which is a version ``up to isomorphism'' of the universal property of the associated sheaf of a (separated) presheaf. (Because of the uniqueness up to strong equivalence one often speaks of \emph{the} associated stack of a prestack $\mathbb{F}$.)

\begin{thm}
For every prestack there exists an associated stack.
\end{thm}
\bew{
Let $\mathbb{F}$ be a prestack. Its associated stack $\hat{\mathbb{F}}$ could be defined briefly as
\[
\hat{\mathbb{F}}(U)=\colimit_{(\mathcal{U})}\mathrm{Des}(\mathcal{U},\mathbb{F})
\]
where $\mathcal{U}$ ranges over covers of $U$. (This colimit should be interpreted with lots of ``up to isomorphisms'', as a pseudo-colimit of categories.) More concretely, define $\hat{\mathbb{F}}(U)$ as follows: an \emph{object} of $\hat{\mathbb{F}}(U)$ consists of
\begin{itemize}
\item a cover $\mathcal{U}=\{U_{\alpha}\}$ of $U$
\item an object $(a,\theta)$ of Des$(\mathcal{U},\mathbb{F})$.
\end{itemize}
For two such objects $(a,\theta)$ in Des$(\mathcal{U},\mathbb{F})$ and $(b,\rho)$ of Des$(\mathcal{V},\mathbb{F})$, for open covers $\mathcal{U}=\{U_{\alpha}\}_{\alpha\in A}$ and $\mathcal{V}=\{V_{\beta}\}_{\beta\in B}$ of $U$, an \emph{arrow} $f:(a,\theta)\rightarrow(b,\rho)$ consists of a system of arrows
\[
f_{\alpha,\beta}:a_{\alpha}|U_{\alpha}\cap V_{\beta}\rightarrow b_{\beta}|U_{\alpha}\cap V_{\beta}
\]
with the property that for any $\alpha,\alpha'\in A$ and $\beta,\beta'\in B$, the diagram
\[
\xymatrix{a_{\alpha}|U_{\alpha\alpha'}\cap V_{\beta\beta'}\ar[r]^{f_{\alpha\beta}} & b_{\beta}|U_{\alpha\alpha'}\cap V_{\beta\beta'} \\ a_{\alpha'}|U_{\alpha\alpha'}\cap V_{\beta\beta'}\ar[u]^{\theta_{\alpha\alpha'}}\ar[r]^{f_{\alpha'\beta'}} & b_{\beta'}|U_{\alpha\alpha'}\cap V_{\beta\beta'}\ar[u]_{\rho_{\beta\beta'}}
}
\]
commutes (where we have omitted the ``$\ldots|U_{\alpha\alpha'}\cap U_{\beta\beta'}$'' for the names of the arrows). The categories $\hat{\mathbb{F}}(U)$ thus defined together have the structure of a fibered category, and there is an evident map $\mathbb{F}\stackrel{w}{\rightarrow}\hat{\mathbb{F}}$ of fibered categories (cf. Remark \ref{rem3} on page \pageref{rem3}). It is now readily checked that $\hat{\mathbb{F}}$ is a stack and $w$ a weak equivalence.
}

\vspace{0.4cm}
\opm{
Recall that any sheaf of categories on $X$ defines a prestack, but this prestack need not be a stack. Moreover, its associated stack is in general no longer a sheaf. (In fact it is hard for a prestack to be a sheaf and a stack at the same time.) On the other hand, every stack is the associated stack of a sheaf of categories. To see this, let $\mathbb{S}$ be a stack. Choose for every open set $U$ a set $R(U)$ of objects in $\mathbb{S}(U)$, so that any object in $\mathbb{S}(U)$ is isomorphic to some object in $R(U)$. Each such open set $U$ defines a sheaf $\subtilde{U}$ on $X$ (whose \'etale space is the inclusion $U\hookrightarrow X$). Now define a sheaf of categories $\mathbb{C}$, consisting of the two sheaves
\[
\emph{Objects}(\mathbb{C})=\coprod_{U\in\mathcal{O}(X)}\coprod_{r\in R(U)}\subtilde{U}
\]
and
\[
\emph{Arrows}(\mathbb{C})=\coprod_{U,V}\coprod\limits_{{r\in R(U),\atop s\in R(V)}} \mathrm{\underline{Hom}}_{\mathbb{S}}(r|U\cap V,s|U\cap V).
\]

There is an obvious map of fibered categories $\mathbb{C}\rightarrow\mathbb{S}$, which is easily seen to be a weak equivalence.
}

\newpage

\section{Gerbes}
We continue to work over a fixed space $X$. A stack $\mathbb{F}$ on $X$ is said to be a \emph{stack of groupoids} if each category $\mathbb{F}(U)$ is a groupoid.

\begin{defi}
A \emph{gerbe} on $X$ is a stack of groupoids $\mathbb{G}$ on $X$ with the following additional properties:
\begin{enumerate}
\item \emph{non-empty:}\quad $X=\bigcup\{U|\mathbb{G}(U)\neq\emptyset\}$\label{ax1gerb}
\item \emph{transitive:}\quad given objects $a,b\in\mathbb{G}(U)$, any point $x\in U$ has a neighborhood $V\subseteq U$ for which there is at least one arrow $a|V\rightarrow b|V$ in $\mathbb{G}(V)$.\label{ax2gerb}
\end{enumerate}
\end{defi}

\noindent{\bf Example:}
The stack $U\mapsto\mathrm{Tor}(\mathcal{G}|U)$ associated to a sheaf of groupoids $\mathcal{G}$ on $X$ (see p. \pageref{asstack}) is a gerbe.

\vspace{0.4cm}
A connected (or transitive) groupoid is characterized up to equivalence by any of its vertex groups. For a \emph{gerbe}, the notion of vertex group is replaced by that of a \emph{band} (``lien'' in {\bf [Gi]}). Before defining the notion of an abstract band, let us describe the band of a gerbe.

Let $\mathbb{G}$ be a gerbe. For an element (object) $a\in\mathbb{G}(U)$ the sheaf $\underline{\mathrm{Hom}}_{\mathbb{G}}(a,a)$ (see Lemma \ref{hom}) is actually a sheaf of groups on $U$, which we denote by $\underline{\mathrm{Aut}}(a)$ (or $\underline{\mathrm{Aut}}_{\mathbb{G}}(a)$ if necessary). Now choose (by axiom \ref{ax1gerb} for gerbes) an open cover $X=\bigcup U_{\alpha}$ and for each $\alpha$ an object $a_{\alpha}\in\mathbb{G}(U_{\alpha})$, thus defining a sheaf of groups $\underline{\mathrm{Aut}}(a_{\alpha})$ on $U_{\alpha}$. Next, by axiom \ref{ax2gerb}, there exists for each $\alpha$ and $\beta$ an open cover $U_{\alpha\beta}=\bigcup U_{\alpha\beta}^{\xi}$ and for each $\xi$ an arrow $f_{\alpha\beta}^{\xi}:a_{\beta}\rightarrow a_{\alpha}$ in $\mathbb{G}(U_{\alpha\beta}^{\xi})$. Conjugation by this arrow defines an isomorphism of sheaves of groups
\[
\lambda_{\alpha\beta}^{\xi}=(f_{\alpha\beta}^{\xi})_*:\underline{\mathrm{Aut}}(a_{\beta})|U_{\alpha\beta}^{\xi}\rightarrow\underline{\mathrm{Aut}}(a_{\alpha})|U_{\alpha\beta}^{\xi}.
\]
This isomorphism $\lambda_{\alpha\beta}^{\xi}$ depends on the choice of $f_{\alpha\beta}^{\xi}$, but two different choices define the same \emph{outer} isomorphism. In particular, on overlaps $U_{\alpha\beta}^{\xi}\cap U_{\alpha\beta}^{\zeta}$, the two isomorphisms $\lambda_{\alpha\beta}^{\xi}$ and $\lambda_{\alpha\beta}^{\zeta}$ define the same outer isomorphism. This means that, for fixed $\alpha$ and $\beta$, the family $\{\lambda_{\alpha\beta}^{\xi}\}_{\xi}$ defines an ``outer isomorphism'' of sheaves on $U_{\alpha\beta}$,
\[
\lambda_{\alpha\beta}:\mathrm{\underline{Aut}}(a_{\beta})|U_{\alpha\beta} \rightarrow\mathrm{\underline{Aut}}(a_{\alpha})|U_{\alpha\beta}
\]
in the sense made precise below, which does not depend on the choice of the $f_{\alpha\beta}^{\xi}$. The system of sheaves of groups $\underline{\mathrm{Aut}}(a_{\alpha})$ and outer isomorphisms $\lambda_{\alpha\beta}$ is called ``the'' \emph{band of the gerbe} $\mathbb{G}$, and is denoted Band$(\mathbb{G})$. (The definition of band, and of isomorphism of bands, will be such that different choices of objects $a_{\alpha}$ define an isomorphic band.)\\

We now introduce the concept of an abstract \emph{band}. First we should make these ``outer isomorphisms'' more precise. Let $K$ and $L$ be two sheaves of groups on $X$. Consider for each $U\subseteq X$ the set Iso$(K|U,L|U)$ of isomorphisms of sheaves of groups on $U$. The group $L(U)$ acts on this set by conjugation, and taking the quotient defines a \emph{presheaf} $U\mapsto\mathrm{Iso}(K|U,L|U)/L(U)$ on $X$. Its associated sheaf is denoted \underline{Out}$(K,L)$ and is referred to as the \emph{sheaf of outer isomorphisms}. A section $\phi\in\underline{\mathrm{Out}}(K,L)(U)$ will be referred to as an outer isomorphism $K|U\rightarrow L|U$. Such a $\phi$ is \emph{represented} by a cover $U=\bigcup U_{\alpha}$ and for each $\alpha$ an isomorphism $\phi_{\alpha}:K|U_{\alpha}\rightarrow L|U_{\alpha}$ (of sheaves of groups on $U_{\alpha}$). These $\phi_{\alpha}$ are compatible in the sense that for each $\alpha$ and $\beta$ there is a cover $U_{\alpha\beta}=\bigcup U_{\alpha\beta}^{\xi}$ such that $\phi_{\alpha}|U_{\alpha\beta}^{\xi}$ and $\phi_{\beta}|U_{\alpha\beta}^{\xi}$ differ by conjugation by an element $\lambda_{\alpha\beta}^{\xi}\in L(U_{\alpha\beta}^{\xi})$.

\vspace{0.4cm}
\opm{
Let $Z(L)$ be the center of $L$. If $H^1(U_{\alpha\beta},Z(L)|U_{\alpha\beta})=0$ for each $U_{\alpha\beta}$, then in the situation above one can choose just one $\lambda_{\alpha\beta}\in L(U_{\alpha\beta})$ such that $\phi_{\alpha}$ and $\phi_{\beta}$ differ on $U_{\alpha\beta}$ by conjugation by $\lambda_{\alpha\beta}$.
}

\vspace{0.4cm}
A \emph{band} $K$ over $X$ is represented by a cover $X=\bigcup U_{\alpha}$, and for each $\alpha$ a sheaf of groups $K_{\alpha}$ on $U_{\alpha}$, together with \emph{outer isomorphisms} $\lambda_{\alpha\beta}:K_{\beta}|U_{\alpha\beta}\rightarrow K_{\alpha}|U_{\alpha\beta}$ satisfying the cocycle condition $\lambda_{\alpha\alpha}=1$ and $\lambda_{\alpha\beta}\lambda_{\beta\gamma}=\lambda_{\alpha\gamma}$ (as outer isomorphisms). The restriction to a finer cover should be viewed as defining the \emph{same} band. If $K=(K_{\alpha},\lambda_{\alpha\beta})$ and $L=(L_{\alpha},\mu_{\alpha\beta})$ are two bands over $X$, represented over the same cover $\{U_{\alpha}\}$, then an isomorphism $K\rightarrow L$ consists of outer isomorphisms $\phi_{\alpha}:K_{\alpha}\rightarrow L_{\alpha}$, compatible on overlaps in the sense that $\phi_{\alpha}\lambda_{\alpha\beta}=\mu_{\alpha\beta}\phi_{\beta}$. This defines the category of bands and isomorphisms over $X$.

\vspace{0.4cm}
\opm{
If each $K_{\alpha}$ is abelian, then the $\lambda_{\alpha\beta}$ satisfy the cocycle condition on the nose, and (because sheaves of groups form a stack!) there is a sheaf of abelian groups $K$ such that $K|U_{\alpha}\simeq K_{\alpha}$ (by an isomorphism compatible with the $\lambda_{\alpha\beta}$). Thus, up to isomorphism, an abelian band is just a sheaf of abelian groups. In particular, this implies that the center $Z(K)$ of a band $K$ is actually (isomorphic to) a sheaf of abelian groups.
}

\vspace{0.4cm}
\opm{
Let $K=(K_{\alpha},\lambda_{\alpha\beta})$ be a general band. Under the simplifying hypothesis that $H^1(U_{\alpha\beta},K_{\alpha}|U_{\alpha\beta})=0$ and $H^2(U_{\alpha\beta},Z(K_{\alpha})|U_{\alpha\beta})=0$, we can represent each $\lambda_{\alpha\beta}$ by an actual isomorphism of sheaves of groups $\lambda_{\alpha\beta}:K_{\beta}|U_{\alpha\beta}\rightarrow K_{\alpha}|U_{\alpha\beta}$. Suppose we have chosen such $\lambda_{\alpha\beta}$. Then similarly, if $H^1(U_{\alpha\beta\gamma},Z(K_{\alpha})|U_{\alpha\beta\gamma})=0$, the cocycle condition $\lambda_{\alpha\beta}\lambda_{\beta\gamma}=\lambda_{\alpha\gamma}$ (initially an equality as outer isomorphisms) will hold in the form of an actual equality of isomorphisms of sheaves on $U_{\alpha\beta\gamma}$, as $\lambda_{\alpha\beta}\lambda_{\beta\gamma}=g_*\lambda_{\alpha\gamma}$, where $g\in K_{\alpha}(U_{\alpha\beta\gamma})$ is a section and $g_*$ is the corresponding inner automorphism. For a nice space $X$ we can choose this cover $\{U_{\alpha}\}$ so that all finite intersections of elements in the cover are either empty or contractible, and these simplifying hypotheses are usually fulfilled.
}

\begin{defi}
Let $K$ be a band on $X$. A \emph{gerbe with band $K$} is a gerbe $\mathbb{G}$ on $X$ together with an isomorphism of bands $\theta:\mathrm{Band}(\mathbb{G})\stackrel{\sim}{\rightarrow} K$. \end{defi}\label{simphyp}

If $(\mathbb{G},\theta)$ and $(\mathbb{G}',\theta')$ are gerbes with band $K$, then a morphism $(\mathbb{G},\theta)\rightarrow(\mathbb{G}',\theta')$ is a morphism of fibered categories $\mu:\mathbb{G}\rightarrow\mathbb{G'}$ for which $\theta'\mu=\theta$ (as maps between bands). Such a morphism is automatically an equivalence. (Recall that any weak equivalence between stacks is strong.) Two such gerbes with band $K$, say $(\mathbb{G},\theta)$ and $(\mathbb{G}',\theta')$, are said to be \emph{equivalent} if there exists a third gerbe $(\mathbb{H},\rho)$ with band $K$ for which there are morphisms
\[
(\mathbb{G},\theta)\leftarrow(\mathbb{H},\rho)\rightarrow(\mathbb{G}',\theta').
\]
The set of equivalence classes of gerbes with band $K$ is denoted 
\[
\mathrm{Gerbes}(X,K).
\]

We will now give a ``cocycle description'' of $\mathrm{Gerbes}(X,K)$ which will lead to a \v{C}ech cohomology (in degree 2) with coefficients in a band. Suppose the band $K$ is given by sheaves of groups $K_{\alpha}$ over $U_{\alpha}$ and outer isomorphisms $\lambda_{\alpha\beta}:K_{\beta}|U_{\alpha\beta}\rightarrow K_{\alpha}|U_{\alpha\beta}$ as before. Let us make the simplifying hypotheses of the previous Remark, so that the $\lambda_{\alpha\beta}$ can each be represented by an \emph{actual} isomorphism (again called $\lambda_{\alpha\beta}$) of sheaves of groups, and let us fix such a choice, where we agree to choose $\lambda_{\alpha\alpha}=\mathrm{identity}$.

Now let $(\mathbb{G},\theta)$ be a gerbe with band $K$. Suppose that the band of $\mathbb{G}$ is represented by objects $a_{\alpha}\in\mathbb{G}(U_{\alpha})$ and arrows $f_{\alpha\beta}:a_{\beta}|U_{\alpha\beta}\rightarrow a_{\alpha}|U_{\alpha\beta}$ in $\mathbb{G}|U_{\alpha\beta}$. (By passing to a common refinement if necessary, we may take the same cover $X=\bigcup U_{\alpha}$ here; and we may choose the arrows $f_{\alpha\beta}$ over $U_{\alpha\beta}$ rather than over a cover $U_{\alpha\beta}^{\xi}$, because of the same simplifying hypothesis -- $H^1(U_{\alpha\beta},\underline{\mathrm{Aut}}(a_{\alpha}))=0$ in this case.) Since $\theta$ is a morphism of bands, for each $U_{\alpha\beta}$ the diagram below

\[
\xymatrix{\underline{\mathrm{Aut}}(a_{\beta}) \ar[d]_{(f_{\alpha\beta})_*} \ar[r]^{\theta_{\beta}} & K_{\beta} \ar[d]^{\lambda_{\alpha\beta}}\\ \underline{\mathrm{Aut}}(a_{\alpha}) \ar[r]^{\theta_{\alpha}} & K_{\alpha}}\label{diag36}
\]

\noindent commutes as a square of outer isomorphisms (where we have omitted ``$\ldots|U_{\alpha\beta}$'' from the notation). By our simplifying hypothesis $(H^1(U_{\alpha\beta},Z\underline{\mathrm{Aut}}(a_{\alpha}))=0)$ we can actually assume that there is one $g\in\underline{\mathrm{Aut}}(a_{\alpha})(U_{\alpha\beta})$ for which $\lambda_{\alpha\beta}\theta_{\beta}=g_*(f_{\alpha\beta})_*\theta_{\alpha}$ (after having represented the $\theta_{\alpha}$'s by actual maps as well). So if we replace $f_{\alpha\beta}$ by $gf_{\alpha\beta}$, we can assume that each square as above commutes on the nose. Also we can assume $f_{\alpha\alpha}=1$ for each $\alpha$. Now write
\[
g_{\alpha\beta\gamma}=\theta_{\alpha}(f_{\alpha\beta}f_{\beta\gamma}f_{\alpha\gamma}^{-1})\in K_{\alpha}(U_{\alpha\beta\gamma}).
\]

\begin{prop}
These $g_{\alpha\beta\gamma}$ satisfy the following cocycle and normality identities,
\begin{enumerate}
\item $\lambda_{\alpha\beta}\lambda_{\beta\gamma}=(g_{\alpha\beta\gamma})_*\lambda_{\alpha\gamma}$ (identity of isomorphisms $K_{\gamma}|U_{\alpha\beta\gamma}\rightarrow K_{\alpha}|U_{\alpha\beta\gamma}$)
\item $g_{\alpha\beta\gamma}g_{\alpha\gamma\delta}=\lambda_{\alpha\beta}(g_{\beta\gamma\delta})g_{\alpha\beta\delta}$ (identity in $K_{\alpha}(U_{\alpha\beta\gamma\delta})$)\label{gabg2}
\item $\lambda_{\alpha\alpha}=1,\ g_{\alpha\alpha\gamma}=1=g_{\alpha\gamma\gamma}$.
\end{enumerate}
\end{prop}

\bew{
\begin{enumerate}
\item By commutativity of the square above,
\begin{eqnarray*}
\lambda_{\alpha\beta}\lambda_{\beta\gamma}&=&\theta_{\alpha}(f_{\alpha\beta})_*(f_{\beta\gamma})_*\theta_{\gamma}^{-1}\\
&=&\theta_{\alpha}(f_{\alpha\beta}f_{\beta\gamma}f_{\alpha\gamma}^{-1})_*(f_{\alpha\gamma})_*\theta_{\gamma}^{-1}\\
&=&(g_{\alpha\beta\gamma})_*\theta_{\alpha}(f_{\alpha\gamma})_*\theta_{\gamma}^{-1}\\
&=&(g_{\alpha\beta\gamma})_*\lambda_{\alpha\gamma}.
\end{eqnarray*}
\item The $f_{\alpha\beta\gamma}:=f_{\alpha\beta}f_{\beta\gamma}f_{\alpha\gamma}^{-1}$ satisfy the usual cocycle condition
\begin{eqnarray*}
f_{\alpha\beta\gamma}f_{\alpha\gamma\delta}&=&f_{\alpha\beta}f_{\beta\gamma}f_{\gamma\delta}f_{\alpha\delta}^{-1}\\
&=&f_{\alpha\beta}(f_{\beta\gamma}f_{\gamma\delta}f_{\beta\delta}^{-1})f_{\alpha\beta}^{-1}f_{\alpha\beta}f_{\beta\delta}f_{\alpha\delta}^{-1}\\
&=&(f_{\alpha\beta})_*(f_{\beta\gamma\delta})f_{\alpha\beta\delta},
\end{eqnarray*}
and this gives \ref{gabg2} after applying $\theta_{\alpha}$.
\end{enumerate}
}

\subsection*{Equivalence of cocycles}
To find the appropriate equivalence relation, let us consider the choices involved in obtaining the 2-cocycle from the gerbe: we have chosen objects $a_{\alpha}$ to represent the band of $\mathbb{G}$, we have chosen \emph{actual} isomorphisms of sheaves $\lambda_{\alpha\beta}$ and $\theta_{\alpha}$ to represent outer isomorphisms and we have chosen $f_{\alpha\beta}:a_{\beta}\rightarrow a_{\alpha}$ in $\mathbb{G}(U_{\alpha\beta})$ which make the diagram on page \pageref{diag36} commute. Let us see how a change in any of these choices affects to ``cocycle'' in the Proposition.\\

First, keeping the $a_{\alpha}$ fixed, we could replace $\lambda_{\alpha\beta}$ by $\mu_{\alpha\beta}$ and keep $\theta_{\alpha}$ the same. Say $\mu_{\alpha\beta}=(k_{\alpha\beta})_*\lambda_{\alpha\beta}$ where $k_{\alpha\beta}\in K_{\alpha}(U_{\alpha\beta})$. Write $e_{\alpha\beta}=\theta_{\alpha}^{-1}(k_{\alpha\beta})$ and replace $f_{\alpha\beta}$ by $\bar{f}_{\alpha\beta}:=e_{\alpha\beta}f_{\alpha\beta}$. Then the diagram

\[
\xymatrix{\underline{\mathrm{Aut}}(a_{\beta}) \ar[d]^{(f_{\alpha\beta})_*} \ar[r]^{\theta_{\beta}} \ar@/_1.6pc/[dd]_{(\bar{f}_{\alpha\beta})_*}& K_{\beta} \ar[d]_{\lambda_{\alpha\beta}} \ar@/^1.6pc/[dd]^{\mu_{\alpha\beta}}\\ \underline{\mathrm{Aut}}(a_{\alpha}) \ar[d]^{(e_{\alpha\beta})_*} \ar[r]^{\theta_{\alpha}} & K_{\alpha} \ar[d]_{(k_{\alpha\beta})_*}\\
\underline{\mathrm{Aut}}(a_{\alpha}) \ar[r]^{\theta_{\alpha}}& K_{\alpha}}
\]
commutes, and the new cocycle
\[
\bar{g}_{\alpha\beta\gamma}:=\theta_{\alpha}(\bar{f}_{\alpha\beta}\bar{f}_{\beta\gamma}\bar{f}_{\alpha\gamma}^{-1})
\]
satisfies
\[
\bar{g}_{\alpha\beta\gamma}k_{\alpha\gamma}=k_{\alpha\beta}\lambda_{\alpha\beta}(k_{\beta\gamma})g_{\alpha\beta\gamma}.
\]

Secondly, we could keep $\lambda_{\alpha\beta}$ the same and replace $\theta_{\alpha}$ by $\rho_{\alpha}$, say
\[
\theta_{\alpha}=\rho_{\alpha}(e_{\alpha})_*\quad\mathrm{where\ }e_{\alpha}\in\underline{\mathrm{Aut}}(a_{\alpha})(U_{\alpha}).
\]
Replace $f_{\alpha\beta}$ by $h_{\alpha\beta}:=e_{\alpha}f_{\alpha\beta}e_{\beta}^{-1}$. Then an easy calculation shows that the new cocycle
\[
\bar{g}_{\alpha\beta\gamma}:=\rho_{\alpha}(h_{\alpha\beta}h_{\beta\gamma}h_{\alpha\gamma}^{-1})
\]
is actually equal to the old one $g_{\alpha\beta\gamma}$.\\

Thirdly, we could keep the $\lambda$'s and the $\theta$'s the same, but replace $f_{\alpha\beta}$ by $z_{\alpha\beta}f_{\alpha\beta}$ where $z_{\alpha\beta}$ is a section of the center $Z(\underline{\mathrm{Aut}}(a_{\alpha}))$ over $U_{\alpha\beta}$. Write $m_{\alpha\beta}=\theta_{\alpha}(z_{\alpha\beta})\in Z(K_{\alpha})$. Then the new cocycle $\bar{g}_{\alpha\beta\gamma}$ is
\begin{eqnarray*}
\bar{g}_{\alpha\beta\gamma}&=&\theta_{\alpha}(z_{\alpha\beta}f_{\alpha\beta}z_{\beta\gamma}f_{\beta\gamma}f_{\alpha\gamma}^{-1}z_{\alpha\gamma}^{-1})\\
&=&m_{\alpha\beta}\lambda_{\alpha\beta}(m_{\beta\gamma})g_{\alpha\beta\gamma}m_{\alpha\gamma}^{-1}.
\end{eqnarray*}

Finally, we could replace $a_{\alpha}$ by $b_{\alpha}$, choose arrows $e_{\alpha}:b_{\alpha}\rightarrow a_{\alpha}$ in $\mathbb{G}(U_{\alpha})$ and replace (our choice for) $\theta_{\alpha}$ by $\bar{\theta}_{\alpha}=\theta_{\alpha}(e_{\alpha})_*$, as in

\[
\begin{array}{cc}
\multicolumn{1}{m{1cm}}{\xymatrix{\underline{\mathrm{Aut}}(b_{\beta})\ar[r]^{(e_{\beta})_*}\ar[d]&\underline{\mathrm{Aut}}(a_{\beta}) \ar[d]_{(f_{\alpha\beta})_*} \ar[r]^{\theta_{\beta}} & K_{\beta} \ar[d]^{\lambda_{\alpha\beta}}\\\underline{\mathrm{Aut}}(b_{\alpha})\ar[r]^{(e_{\alpha})_*} & \underline{\mathrm{Aut}}(a_{\alpha}) \ar[r]^{\theta_{\alpha}} & K_{\alpha}}} &
\multicolumn{1}{m{3cm}}{(over $U_{\alpha\beta})$}
\end{array}
\]

Then (with the same $\lambda$'s) each $f_{\alpha\beta}$ is replaced by $\bar{f}_{\alpha\beta}=e_{\alpha}^{-1}f_{\alpha\beta}e_{\beta}$. However, the new cocycle $\bar{g}_{\alpha\beta\gamma}:=\bar{\theta}_{\alpha}(\bar{f}_{\alpha\beta}\bar{f}_{\beta\gamma}\bar{f}_{\alpha\gamma}^{-1})$ hasn't changed and equals the old one $g_{\alpha\beta\gamma}$. So, summing up:

\begin{defi}
Let $K$ be a band on $X$, and let $\mathcal{U}=\{U_{\alpha}\}$ be an open cover of $X$: A \emph{(normal) 2-cocycle} over $\mathcal{U}$ with values in $K$ consists of 
\begin{itemize}
\item a representation of $K$ by sheaves $K_{\alpha}$ on $U_{\alpha}$ and actual (not outer) isomorphisms $K_{\beta}|U_{\alpha\beta}\stackrel{\lambda_{\alpha\beta}}{\rightarrow}K_{\alpha}|U_{\alpha\beta}$.
\item elements $g_{\alpha\beta\gamma}\in K_{\alpha}(U_{\alpha\beta\gamma})$ such that the cocycle conditions of the Proposition all hold.
\end{itemize}
\end{defi}

Two such cocycles $(\lambda_{\alpha\beta},g_{\alpha\beta\gamma})$ and $(\mu_{\alpha\beta},h_{\alpha\beta\gamma})$ are \emph{equivalent} if there are $k_{\alpha\beta}\in K_{\alpha}(U_{\alpha\beta})$ such that
\begin{itemize}
\item $\mu_{\alpha\beta}=(k_{\alpha\beta})_*\lambda_{\alpha\beta}$
\item $h_{\alpha\beta\gamma}=k_{\alpha\beta}\lambda_{\alpha\beta}(k_{\beta\gamma})g_{\alpha\beta\gamma}k_{\alpha\gamma}^{-1}$.
\end{itemize}
The set of equivalence classes is denoted
\[
\check{H}^2(\mathcal{U},K)
\]
and the colimit over open covers
\[
\check{H}^2(X,K)=\colimit_{\mathcal{U}}\check{H}^2(\mathcal{U},K).
\]
Elements are referred to as non-abelian \v{C}ech cohomology classes.

\begin{thm}
Let $K$ be a band on $X$. There is a bijective correspondence
\[
\mathrm{Gerbes}(X,K)\cong\check{H}^2(X,K)
\]
between equivalence classes of gerbes with band $K$ and non-abelian cohomology classes.
\end{thm}

\bew{
From left to right, we have already explained how to obtain an equivalence class of cocycles from a gerbe. From right to left, one can obtain an inverse construction as follows. For a cocycle $\{\lambda_{\alpha\beta},g_{\alpha\beta\gamma}\}$ one can construct a sheaf of groupoids $\mathcal{G}$ with non-empty connected stalks. The associated stack of this sheaf provides the required inverse construction. The sheaf of objects $\mathcal{G}_0$ of $\mathcal{G}$ is the sum
\[
\mathcal{G}_0=\coprod_{\alpha}\subtilde{U}_{\alpha}
\]
where $\subtilde{U}_{\alpha}$ is the sheaf corresponding to the \'etale space $U_{\alpha}\hookrightarrow X$. The sheaf of arrows $\mathcal{G}_1$ of $\mathcal{G}$ is the sum
\[
\mathcal{G}_1=\coprod_{\alpha,\beta}K_{\alpha}|U_{\alpha\beta}.
\]
In other words, an object of the stalk $\mathcal{G}_x$ at $x$ is a pair $(x,\alpha)$ where $x\in U_{\alpha}$, and an arrow $(x,\beta)\stackrel{k}{\rightarrow}(x,\alpha)$ is a point $k$ in the stalk of $K_{\alpha}$ at $x$. Composition of two such, $(x,\beta)\stackrel{k}{\rightarrow}(x,\alpha)$ and $(x,\gamma)\stackrel{l}{\rightarrow}(x,\beta)$, is defined to be the arrow
\[
k\cdot\lambda_{\alpha\beta}(l)\cdot g_{\alpha\beta\gamma}(x):(x,\gamma)\rightarrow(x,\alpha),
\]
where the product is taken in the group $(K_{\alpha})_x$. The cocycle relation implies that this is an associative operation. It is now straightforward to check that the two constructions are mutually inverse.
}

\vspace{0.4cm}
\opm{
Suppose that the band $K$ is a \emph{sheaf} of groups, so that we can take $K_{\alpha}=K|U_{\alpha}$ and each $\lambda_{\alpha\beta}$ the identity. Then for a cocycle
$\{\lambda_{\alpha\beta},g_{\alpha\beta\gamma}\}$ with these $\lambda$'s, each $g_{\alpha\beta\gamma}$ lies in the center of $K_{\alpha\beta\gamma}$ and the usual abelian cocycle relation
\[
g_{\alpha\beta\gamma}g_{\alpha\gamma\delta}=g_{\beta\gamma\delta}g_{\alpha\beta\delta}
\]
holds. So gerbes whose band is isomorphic to a sheaf of groups correspond to usual (abelian) 2-cocycles with coefficients in the center $Z(K)$.
}

\vspace{0.4cm}
\opm{
Suppose we are given a band $K$ on $X$, represented by an open cover $\mathcal{U}=\{U_{\alpha}\}$ of $X$, sheaves $K_{\alpha}$ on $U_{\alpha}$, and isomorphisms $\lambda_{\alpha\beta}$ as before. One may ask whether there exists a gerbe with band $K$. To answer this question, let us again make the simplifying hypotheses that $H^1(U_{\alpha\beta},K|U_{\alpha\beta})$ and $H^2(U_{\alpha\beta},Z(K_{\alpha})|U_{\alpha\beta})$ are trivial (as in the Remark preceding Definition \ref{simphyp}), so that the $\lambda_{\alpha\beta}$ can be represented by \emph{actual} isomorphisms $\lambda_{\alpha\beta}:K_{\beta}|U_{\alpha\beta}\rightarrow K_{\alpha}|U_{\alpha\beta}$ of sheaves of groups, and that the cocycle condition for the $\lambda_{\alpha\beta}$ can be represented in the form $\lambda_{\alpha\beta}\lambda_{\beta\gamma}=(g_{\alpha\beta\gamma})_*\lambda_{\alpha\gamma}$ for sections $g_{\alpha\beta\gamma}\in K_{\alpha}(U_{\alpha\beta\gamma})$. Then
\[
(g_{\alpha\beta\gamma}g_{\alpha\gamma\delta})_*=(\lambda_{\alpha\beta}(g_{\beta\gamma\delta})g_{\alpha\beta\delta})_*
\]
so the ``boundary'' of the non-abelian cochain $\{g_{\alpha\beta\gamma}\}$,
\[
\xi_{\alpha\beta\gamma\delta}=g_{\alpha\beta\gamma}g_{\alpha\gamma\delta}\lambda_{\alpha\beta}(g^{-1}_{\alpha\beta\delta})g^{-1}_{\beta\gamma\delta}
\]
is a section of $Z(K_{\alpha})$, and we obtain a class $[\xi]$ where $\xi=\{\xi_{\alpha\beta\gamma\delta}\}$, in the (usual, abelian) \v{C}ech cohomoloy $H^3(\mathcal{U},Z(K))$. If this cohomology is trivial, then $\xi=d\zeta$ for some $\zeta=\{\zeta_{\alpha\beta\gamma}\}$, and we can replace $g_{\alpha\beta\gamma}$ by $h_{\alpha\beta\gamma}:=\zeta_{\alpha\beta\gamma}g_{\alpha\beta\gamma}$, to obtain an actual cocycle $h$, i.e.
\[
h_{\alpha\beta\gamma}h_{\alpha\gamma\delta}=\lambda_{\alpha\beta}(h_{\beta\gamma\delta})h_{\alpha\beta\delta}.
\]
Thus, in conclusion, if $K$ is a band represented over a cover $\mathcal{U}$, and moreover $H^3(\mathcal{U},Z(K))=0$, then there exists a gerbe with band $K$.
}  

\newpage

\section{Bundle gerbes}
We now work in the context of $C^{\infty}$-manifolds. Recall that a \emph{groupoid} (in Sets) is a small category in which each arrow is an isomorphism. Let us write $G$ for the set of arrows of a groupoid and $M$ for the set of objects. Then the groupoid can also be described in terms of its structure maps

\begin{equation}
\xymatrix{G\times_M G \ar[r]^{\qquad m} & G\ar[r]^{i}&G\ar@<0.5ex>[r]^{s}\ar@<-0.5ex>[r]_{t}&M\ar[r]^{u}&G
}\label{struct}
\end{equation}

\noindent (here $u(x)=1_x$ is the unit at $x\in M,\ i(g)=g^{-1}$ is the inverse of $g\in G$, $s$ and $t$ denote the source and target and $m(g,h)=gh=g\circ h$ is the composition of two arrows $g$ and $h$ with $s(g)=t(h)$.) This groupoid is called a \emph{Lie groupoid (over $M$)} if $M$ and $G$ are smooth manifolds, $s$ and $t$ are submersions and all other structure maps are smooth. In the sequel we write $g:x\rightarrow y$ to denote an element $g\in G$ with $s(g)=x$ and $t(g)=y$.

\vspace{0.4cm}
\vb{
\item If $M$ is a point, this definition reduces to that of a Lie group. More generally, a Lie groupoid  with $s=t$ is a \emph{bundle of Lie groups}. Conversely, in an arbitrary Lie groupoid $G$ over $M$, the set of arrows $x\rightarrow x$ from a given object to itself, has the structure of a Lie group, denoted $G_x$ or $\mathrm{Aut}_G(x)$. It is the \emph{isotropy group} at $x$ (or the \emph{vertex group} at $x$).
\item If $\pi:M\rightarrow X$ is a submersion, the fibered product $G=M\times_X M$ is a Lie groupoid, with exactly one arrow $(x,y):x\rightarrow y$ iff $\pi(x)=\pi(y)$ (for $x,y\in M$). This applies in particular to the case where $U=\coprod U_{\alpha}$ for an open cover $X=\bigcup U_{\alpha}$, with $\pi:U\rightarrow X$ the obvious map.\label{exlie2}
\item Let $X$ be a ``base'' manifold. A \emph{family of groupoids} over $X$ is a Lie groupoid 
($\xymatrix{G \ar@<0.5ex>[r]^{s} \ar@<-0.5ex>[r]_{t}&M}$, etc), equipped with a submersion $\pi:M\rightarrow X$ for which $\pi s=\pi t$. For such a family, each fiber over a point $x\in X$ is a Lie groupoid, denoted $\xymatrix{G_x \ar@<0.5ex>[r] \ar@<-0.5ex>[r]&M_x}$. A \emph{sheaf of groupoids} on $X$ is the same thing as a family with the additional property that the maps $\pi:M\rightarrow X$ and $\pi s=\pi t:G\rightarrow X$ are \'etale (local diffeomorphisms). This implies that all the structure maps in (\ref{struct}) are \'etale as well. (A groupoid with the latter property is called an \emph{\'etale groupoid}.)
}

A \emph{morphism of Lie groupoids} is a smooth functor. In this section, we will in particular be interested in morphisms $G\stackrel{\phi}{\rightarrow}H$ between Lie groupoids over the same manifold $M$ of objects. (Such a morphism is given by a smooth map $\phi:G\rightarrow H$ with $s\phi=s,t\phi=t,\ \phi u=u,\ \phi i=i$, and $m(\phi\times\phi)=\phi m$; cf. (\ref{struct}).)

An \emph{extension of Lie groupoids} over $M$ is a sequence of groupoids and morphisms over $M$,
\[
K\stackrel{j}{\rightarrow}G\stackrel{\phi}{\rightarrow}H,
\]
where $\phi$ is a surjective submersion and $j$ is an embedding onto the submanifold $\mathrm{Ker}(\phi)=\{g\in G|\phi(g)\mathrm{\ is\ a\ unit}\}$. Note that $s(g)=t(g)$ if $g\in\mathrm{Ker}(\phi)$, so that $K$ is in fact a bundle of Lie groups over $M$.

Now let $L$ be a fixed Lie group. For a manifold $M$, we denote by $\underline{L}$ (or $\underline{L}_M$ if necessary) the trivial bundle $L\times M\rightarrow M$ of Lie groups over $M$.

\begin{defi}
An \emph{$L$-bundle gerbe} over a manifold $X$ is given by
\begin{itemize}
\item a surjective submersion $\pi:M\twoheadrightarrow X$
\item an extension $\beta=(\underline{L}\stackrel{j}{\rightarrow}G\stackrel{\phi}{\twoheadrightarrow}M\times_X M)$ over $M$. (Here $M\times_X M$ is the groupoid defined by $\pi$, cf. Example \ref{exlie2}.)
\end{itemize}
\end{defi}
An equivalent definition of an $L$-bundle gerbe is a family of groupoids over $X$,
\[
\xymatrix{G \ar@<0.5ex>[r]^{s} \ar@<-0.5ex>[r]_{t}&M \ar[r]^{\pi}&X}
\]
with the additional properties that
\begin{enumerate}
\item $(s,t):G\rightarrow M\times_X M$ is a surjective submersion;\label{addprop1}
\item there is an isomorphism of Lie groupoids
\[
j_m:L\rightarrow \mathrm{Aut}_G(m)
\]
which identifies each vertex group in $G$ with $L$ and moreover, this isomorphism is required to vary smoothly in $m$.\label{addprop2}
\end{enumerate}
Note that \ref{addprop1} expresses in a strong sense that each fiber $G_x$ is non-empty and connected; so \ref{addprop1} and \ref{addprop2} combined sound like a variation on the definition of ``gerbe with band $L$''. Here is a way to make the connection precise.

Suppose we are given an $L$-bundle gerbe as above. Define a sheaf of groupoids $\mathcal{G}$ over $X$ as follows. Let $\mathcal{G}_0$ be the sheaf of sections of $\pi:M\rightarrow X$. The sheaf $\mathcal{G}_1$ of arrows is defined as follows. For two such sections $a,b\in \mathcal{G}_0(U)$ over an open $U\subseteq X$, an arrow $g:a\rightarrow b$ in $\mathcal{G}_1(U)$ is a section $g:U\rightarrow G$ with $sg=a$ and $tg=b$. Note that the stalk $\mathcal{G}_x$ of $\mathcal{G}$ at $x\in X$ is non-empty (because $\pi$ is a surjective submersion) and connected (because $(s,t)$ is a surjective submersion onto $M\times_X M$). Thus, the associated stack $\hat{\mathcal{G}}$ of $\mathcal{G}$ is a gerbe. What is its band?

The band of $\hat{\mathcal{G}}$ can be described directly in terms of $\mathcal{G}$. To this end, let $X=\bigcup U_{\alpha}$ be an open cover of $X$, for which there exist sections $a_{\alpha}:U_{\alpha}\rightarrow M$ of $\pi$. Since $G\rightarrow M\times_X M$ is a surjective submersion, each point $x\in U_{\alpha\beta}$ has a neighborhood $V$ for which there exists an arrow $a_{\beta}|V\rightarrow a_{\alpha}|V$ in $\mathcal{G}(V)$. If we choose the cover such that each non-empty intersection $U_{\alpha\beta}$ is contractible, then $H^1(U_{\alpha\beta},\subtilde{L})$ is trivial and we can in fact find a section $g_{\alpha\beta}\in \mathcal{G}(U_{\alpha\beta})$,
\[
g_{\alpha\beta}:a_{\beta}|U_{\alpha\beta}\rightarrow a_{\alpha}|U_{\alpha\beta},\quad\mathrm{in\ }\mathcal{G}(U_{\alpha\beta}).
\]
The (inverse of the) map $j$ induces an isomorphism of sheaves of groups on $U_{\alpha}$, denoted
\[
\theta_{\alpha}:\underline{\mathrm{Aut}}(a_{\alpha})\rightarrow\subtilde{L}|U_{\alpha},
\]
where $\subtilde{L}$ is the sheaf of sections of the trivial bundle $\underline{L}_X$. The isomorphisms $\lambda_{\alpha\beta}$ are now defined by the requirement that
\[
\xymatrix{{\underline{\mathrm{Aut}}(a_{\alpha})|U_{\alpha\beta}} \ar[r]^{\theta_{\alpha}}&{\subtilde{L}}|U_{\alpha\beta}\\
{\underline{\mathrm{Aut}}(a_{\beta})|U_{\alpha\beta}}\ar[u]^{(g_{\alpha\beta})_*} \ar[r]^{\theta_{\beta}}&{\subtilde{L}}|U_{\alpha\beta}\ar[u]_{\lambda_{\alpha\beta}}
}
\]
commutes. Then $\{\subtilde{L}|U_{\alpha},\lambda_{\alpha\beta}\}$ is the band of $\mathcal{G}$ (or of $\hat{\mathcal{G}}$). It is some twisted form of the sheaf of groups $\subtilde{L}$. We denote it by $\subtilde{L}^{\beta}$ where $\beta=(\underline{L}\rightarrow G\rightarrow M\times_X M)$ is the extension.

\vspace{0.4cm}
\opm{Suppose $L$ is abelian. Then the extension
\[
L\stackrel{j}{\rightarrow}G\stackrel{\phi}{\rightarrow}M\times_X M
\]
induces a (left, say) action of $M\times_X M$ on $L$ in the obvious way: any pair $(m,n)\in M\times_X M$ gives an isomorphism $(m,n)_*:L\rightarrow L$, defined by $j((m,n)_*(l))=gj(l)g^{-1}$ where $g:m\rightarrow n$ is any arrow in $G$; this is independent of $g$ if $L$ is abelian. (As usual we call the extension \emph{central} if this action is trivial, i.e. each $(m,n)_*$ is the identity.) Now for the $g_{\alpha\beta}:a_{\beta}\rightarrow a_{\alpha}$ as above, this action induces in particular a map $(a_{\beta},a_{\alpha})_*:\subtilde{L}|U_{\alpha\beta}\rightarrow\subtilde{L}|U_{\alpha\beta}$ and this is the $\lambda_{\alpha\beta}$ above. Thus, \emph{the band structure} given by the $\lambda_{\alpha\beta}$ in the non-abelian case \emph{replaces the action} of $M\times_X M$ on $L$ in the abelian case.
}

\vspace{0.4cm}
\opm{
The gerbe $\hat{\mathcal{G}}$ with band $\subtilde{L}^{\beta}$ is characterized by a class
\[
[\beta]\in\check{H}^2(\mathcal{U},\subtilde{L}^{\beta}),
\]
as explained in the previous lecture, namely the class of the cocycle \linebreak $\theta_{\alpha}(g_{\alpha\beta}g_{\beta\gamma}g_{\alpha\gamma}^{-1})$ in the present notation. By the construction at the end of the previous lecture, it is possible to reconstruct, from this cocycle, the restriction of the extension to $U=\coprod U_{\alpha}$, as in the pullback
\[
\xymatrix{\underline{L}_U \ar@{^{(}->}[d] \ar[r] & G_U \ar@{^{(}->}[d] \ar[r] & U\times_X U \ar@{^{(}->}[d] \\
\underline{L}_M \ar[r] & G \ar[r] & M\times_X M.
}
\]
}\label{p45}

\noindent However, for a general extension $\underline{L}_U\rightarrow H\rightarrow U\times_X U$ of groupoids over $U$, there is an obstruction in $H^2(M,\subtilde{L})$ to being the pullback of an extension $\underline{L}_M\rightarrow G\rightarrow M\times_X M$ of groupoids over $M$ (see {\bf [Mo]}).

\vspace{0.4cm}
\noindent{\bf Example:} Let $L\rightarrow G\rightarrow H$ be an extension of Lie groups, and let $P\rightarrow X$ be a principal $H$-bundle. Then each of $H,G$ and $L$ act on $P$ and give rise to an extension of translation groupoids
\[
L\ltimes P\rightarrow G\ltimes P\rightarrow H\ltimes P
\]
over $P$. (For example, arrows $p\rightarrow q$ in $G\ltimes P$ are $g\in G$ with $gp=q$.) Since the $H$-action is principal and the $L$-action is trivial, this sequence is isomorphic to
\[
\beta:\underline{L}\rightarrow G\times P\rightarrow P\times_X P,
\]
so we have an $L$-bundle gerbe over $X$ and hence a class in $H^2(\mathcal{U},\subtilde{L}^{\beta})$. This class vanishes if $P$ is induced from a principal $G$-bundle. The twisting $\subtilde{L}^{\beta}$ of $L$ is related to that of the original extension and in particular $\beta$ is central if $L\rightarrow G\rightarrow H$ is.\\

In the literature, the following special class of gerbes is often discussed under the name ``bundle gerbe''.
\begin{defi}
A \emph{circle bundle gerbe} over $X$ is given by a surjective submersion $\pi:M\twoheadrightarrow X$ and a \emph{central} extension
\[
\beta:\underline{S}^1\stackrel{j}{\rightarrow}G\stackrel{\phi}{\rightarrow}M\times_X M.
\]
\end{defi}
For such a circle bundle gerbe $\beta$, one obtains a characteristic class $[\beta]\in\check{H}^2(\mathcal{U},\subtilde{S}^1)$. If the cover $\mathcal{U}$ of $X$ is well chosen (contractible intersections $U_{\alpha_1}\cap\ldots\cap U_{\alpha_n}$), then the \v{C}ech cohomology agrees with sheaf cohomology and the long exact cohomology sequence of $0\rightarrow\underline{\mathbb{Z}}\rightarrow\subtilde{\mathbb{R}}\rightarrow\subtilde{S}^1\rightarrow 0$ shows that $\check{H}^2(\mathcal{U},\subtilde{S}^1)=H^2(X,\subtilde{S}^1)=H^3(X,\mathbb{Z})$ so the bundle gerbe has its class in $H^3(X,\mathbb{Z})$.

\vspace{0.4cm}
To conclude this course we give, following {\bf [Mu]}, a construction of this class in De Rham cohomology. To this end, first recall the notion of a connection on a principal $S^1$-bundle $P\rightarrow Y$, which can be defined either in terms of $S^1$-equivariant ``horizontal lifts'' $P\times_Y TY\rightarrow TP$ or in terms of connection 1-forms $\omega$. Recall that a connection $\omega$ on $P\rightarrow Y$ pulls back along a map $f:Z\rightarrow Y$ to a connection on $f^*(P)$, with connection 1-form $f^*(\omega)$ (we also write $f$ for the map $f^*(P)\rightarrow P$). Recall also that for two principal bundles $P\rightarrow Y$ and $Q\rightarrow Y$, their contracted product $P\otimes_{S^1}Q$, obtained from the fibered product $P\times_Y Q$ by identifying $(p\theta,q)$ with $(p,q\theta)$ for $\theta\in S^1$, is again a principal bundle. Connections on $P$ and $Q$ induce a connection on $P\otimes_{S^1}Q$ by taking the evident horizontal lift in each component. If $\omega$ and $\eta$ are the connection 1-forms, then the 1-form for this induced connection is denoted $\omega + \eta$ (it sends a tangent vector $v\otimes w\ (\mathrm{where\ }v\in T_p(P),w\in T_q(Q))$ to $\omega_p(v)+\eta_q(w)$).\\

For the rest of this section, let us fix a circle bundle gerbe
\[
\beta:\underline{S}^1\stackrel{j}{\rightarrow}G\stackrel{(s,t)}{\rightarrow}M\times_X M
\]
over $X$. Observe that $G\stackrel{(s,t)}{\rightarrow}M\times_X M$ is a principal bundle over $M$. The left and right actions agree,
\[
j_y(\theta)g=gj_x(\theta)\quad(\mathrm{for\ }x\stackrel{g}{\rightarrow}y,\theta\in S^1)
\]
because the extension is central. We will write $g\theta$ for $gj_x(\theta)$. Note that the composition of the groupoid
\[
G\times_M G\stackrel{\gamma}{\rightarrow}G\quad(g,h)\mapsto hg
\]
can be viewed as an \emph{isomorphism} of principal bundles over $M\times_X M\times_X M$,
\begin{equation}
\pi_{12}^*(G)\otimes\pi_{23}^*(G)\stackrel{\gamma}{\rightarrow}\pi_{13}^*(G)
\end{equation}\label{pieq}

\noindent with fiber over $(x,y,z)$ mapping $g\otimes h$ to $hg$ where $x\stackrel{g}{\rightarrow}y$ and $y\stackrel{h}{\rightarrow}z$.

\begin{defi}
A \emph{(groupoid) connection} on the bundle gerbe $\beta$ as above is a connection on the principal $S^1$-bundle $G\stackrel{(s,t)}{\rightarrow}M\times_X M$, which preserves the groupoid structure: the composition and inverse of horizontal lifts is again horizontal.
\end{defi}

We will show in a moment that such connections always exist. For now, suppose we are given such a groupoid connection in terms of a 1-form $\omega$. Then the map $\gamma$ in (\ref{pieq}) will be an isomorphism of bundles with connection, respectively $\pi_{12}^*(\omega)+\pi_{23}^*(\omega)$ and $\pi_{13}^*(\omega)$. Let $\kappa$ be the curvature; it is the closed 2-form on $M\times_X M$ defined by $(s,t)^*(\kappa)=d\omega$. Since $\gamma$ is an isomorphism, the two bundles in (\ref{pieq}) have the same curvature, so
\begin{equation}
\pi_{12}^*(\kappa)+\pi_{23}^*(\kappa)=\pi_{13}^*(\kappa).
\end{equation}\label{pieq2}

Now choose a cover $M=\bigcup U_{\alpha}$, let $U=\coprod U_{\alpha}$, and pull everything back along $U\rightarrow M$ (as in the diagram on page \pageref{p45}). Then (the restriction of) $\kappa$ is a closed 2-form on $U\times_X U$, still satisfying the cocycle condition (\ref{pieq2}), on $U\times_X U\times_X U$. So by Mayer-Vietoris (see {\bf [BT]}), there is a 2-form $\lambda$ on $U$ for which $\kappa=\pi_2^*(\lambda)-\pi_1^*(\lambda)$, where $\pi_1,\pi_2:U\times_X U\rightarrow U$ are the projections. Since $d\kappa=0$ we have $\pi_1^*(d\lambda)=\pi_2^*(d\lambda)$, so by Mayer-Vietoris again, $d\lambda=\pi^*(\xi)$ for a closed 3-form $\xi$ on $X$ (where $\pi:U\rightarrow X$). Then $[\xi]\in H^3_{DR}(X)$ and (by the usual tic-tac-toe argument in the \v{C}ech-De Rham complex) it can be shown that $[\xi]$ is integral and coincides with the class constructed earlier.\\

Finally, let us prove that any bundle gerbe $\beta$ as above carries at least one groupoid connection. The problem is local in $X$, by the usual partition of unity argument, so we may assume $\pi:M\rightarrow X$ has a section $X\stackrel{a}{\rightarrow}M$. Let $P_a$ be the pullback of $G\stackrel{(s,t)}{\rightarrow}M\times_X M$ along $(a\pi,\mathrm{id}):M\rightarrow M\times_X M$, so
\[
(P_a)_x=\{\mathrm{arrows\ } a\pi(x)\rightarrow x \},\quad x\in M.
\]
Then $P_a$ is a principal $S^1$-bundle over $M$. Let $\eta$ be an ordinary connection on $P_a$, providing for each $g:a\pi(x)\rightarrow x$ and each $v\in T_x(M)$ a horizontal lift $\tilde{v}_g$. Let $P_a^{\vee}$ be the same space $P_a$ with opposite right $S^1$ action, $(g,\theta)\mapsto g\theta^{-1}$, so that $-\eta$ is a connection on $P_a^{\vee}$ with the \emph{same} horizontal lifts. Now consider the isomorphism of principal bundles
\[
\pi_1^*(P_a^{\vee})\otimes\pi_2^*(P_a)\stackrel{\delta}{\rightarrow}G
\]
defined for $g:a\pi(x)\rightarrow x$ and $h:a\pi(y)\rightarrow y$ by
\[
\delta(g\otimes h)=hg^{-1}.
\]
Then $\delta_*(\pi_1^*(-\eta)+\pi_2^*(\eta))$ is a groupoid connection on $G$. Indeed, for $k:x\rightarrow y$ and tangent vectors $(v,w)\in T_{x,y}(M\times_X M)$, factor $k$ as $x\stackrel{g^{-1}}{\rightarrow}a\pi(x)=a\pi(y)\stackrel{h}{\rightarrow}y$; then the horizontal lift of $(v,w)$ at $k$ is $\tilde{w}_h\circ\tilde{v}_g^{-1}$ (from this formula it is evident that it is a groupoid connection).

\newpage

{\normalfont\LARGE\bfseries References}

\vspace{0.4cm}
\begin{enumerate}
\item[{\bf [BT]}] R. Bott and L.W. Tu, Differential Forms in Algebraic Topology, Springer-Verlag, 1982.\\
\item[{\bf [Br]}] L. Breen, On the classification of 2-gerbes and 2-stacks, Ast\'erisque 225, (1994).\\
\item[{\bf [Gi]}] J. Giraud, Cohomologie non ab\'elienne, Springer-Verlag, 1971.\\
\item[{\bf [Mo]}] I. Moerdijk, On the classification of regular Lie groupoids, preprint math.DG/0203099.\\
\item[{\bf [Mu]}] M. Murray, Bundle gerbes, J. London Math. Soc. 54 (1996).
\end{enumerate}

\end{document}